%Feb 24, 2209

\magnification=1200
\centerline{\bf Pluri-polarity in almost complex structures.}
\vskip1truein
\centerline{Jean-Pierre Rosay}
\vfill
\noindent
Abstract: $J$-holomorphic curves are $-\infty$ sets of $J$-plurisubharmonic functions,
with a singularity of LogLog type, but it is shown that in general they are not
$-\infty$ sets of $J$-plurisubharmonic functions with Logarithmic singularity
(i.e. non-zero Lelong number).
Some few additional remarks on pluripolarity in almost complex structures are made.
\vskip2truein\noindent
A.M.S. Classification: 32Q60, 32Q65, 32U05
\eject

\noindent
INTRODUCTION. $J$-plurisubharmonic functions with
poles (at which the function is $-\infty$) have already played a role in almost complex 
analysis. They  has been used in the study of the Kobayashi metric ([G-S] [I-R]) for getting
an efficient control of $J$-holomorphic discs. In a work in progress with S. Ivashkovich,
applications to uniqueness problems are given. The first pluripolarity result is due to Chirka
who showed that if $J$ is a ${\cal C}^1$ almost complex structure defined near 0 in ${\bf C}^n$
and if $J(0)=J_{st}$ (the standard complex structure), then for $A>0$ large enough 
$\log |z|~+A|z|$ is $J$-plurisubharmonic near 0. A complete proof has been written in [I-R] (Lemma 1.4
page 2400). The function
$-\log |\log |z||$ is also $J$-plurisubharmonic near 0. Although a `log-log singularity'
(zero Lelong number!) is much less interesting and has less applications 
that a `log singularity', functions with
log-log singularity were introduced in [Ro1] to show pluripolarity of $J$-holomorphic curves.
Later, Elkhadhra [E] generalized the result to show pluri-polarity of $J$-holomorphic submanifolds,
again with a `loglog' singularity of the function.
\bigskip\noindent
In [Ro1] there has been an error in stating the smoothness hypotheses, all smoothness
requirements in the statements have to be increased by 1 (${\cal C}^{k+1}$ instead of
${\cal C}^k$). Indeed on line 11-  page 663 it is claimed that because 
$J[Y',JY'](Z',0)=0$, one has $|J[Y',JY'] (Z',Z'')|\leq C |Z''||Y'|^2$.
This is correct if $J$ is of class ${\cal C}^2$ (or at least ${\cal C}^{1,1}$), but ${\cal C}^1$
smoothness, hence continuity of $[Y',JY']$, is not enough.
\bigskip\noindent The really original part of the paper is an example in dimension 4 showing
that {\it in general, $J$-holomorphic curves are not $-\infty$ set of any $J$-plurisubharmonic
function with a logarithmic singularity}. The precise statement is given in Proposition 2.
\bigskip\noindent 
The plan of the paper is as follows:
\medskip\noindent
Part 0 is a naive discussion, in one complex variable, of the functions $\log |z|$ and
$\log |\log |z||$, and of perturbations of the Laplacian. Since the theory of almost complex
structures is in part the theory of small perturbations of the standard complex structure
(and thus of the Laplacian), it seems to me that this preliminary discussion leads to a
better
understanding of the problems that are encountered.
\medskip\noindent
In Part I, we develop the calculus tools to be used in part II and III, for studying $J$-plurisubharmonicity.
\medskip\noindent
Part II is a proof of the pluripolarity of $J$-holomorphic curves established in [Ro1] (with the smoothness
hypothesis
this time correctly stated!). Although the result is not new, the proof is written very differently,
in a way such that I think that it is much easier to see  how things happen. For me, it has been a 
needed preliminary before
getting the counterexample that is given in Part III. 
\medskip\noindent
Part III itself consists of two parts. In the first one one shows the non existence of functions
with log singularity under natural reasonable hypotheses on the behavior of the function. Then in the second part
the non existence result is proved in full generality. For that, we shall need 3 Lemmas in measure
theory whose proof is given in an Appendix.
\medskip\noindent
Part IV treats the case when there is better matching with the standard complex structure,
when the Nijenhuis tensor vanishes along a $J$-holomorphic disc. Things become much easier
and functions with logarithmic poles exist.
\medskip\noindent
Part V is devoted to the very simple remark that in dimension 4 a 2 dimensional surface cannot
be pluripolar if it is not a $J$-holomorphic curve. In higher dimensions $2n$, the corresponding simple
remark is for generic n-dimensional manifolds. Our elementary argument avoids any construction of
'Pinchuck discs' with part of their boundary attached to the submanifold.
\bigskip\noindent
We shall take the following definition of $J$-plurisubharmonicity:
\medskip\noindent
{\bf Definition.} {\it Let $J$ be an almost complex structure defined
on some smooth manifold, of class at least ${\cal C}^{1,\alpha}$ for some
$0<\alpha <1$.
An upper semi-continuous (hence continuous)
function $\lambda$, defined on that set, with values in $[-\infty , +\infty )$ will be said to
be $J$-plurisubharmonic if $\lambda $ is of class ${\cal C}^2$ off its $-\infty$ set and if its
restriction to any $J$-holomorphic curve is subharmonic}.
\medskip\noindent
The hypothesis that $J$ is of class ${\cal C}^{1,\alpha}$ insures that all the $J$-holomorphic discs
are of class ${\cal C}^2$
There is an easy characterization of $J$-plurisubharmonicity in terms of some kind
of Levi form (see e.g. [I-R], Corollary 1.1 page 2400), that we shall not use here,
except briefly (and just for simplification) in parts IV and V. The theory for non ${\cal C}^2$
functions is still not completely satisfactory, although there is the interesting work of Pali
[Pa].
\bigskip\noindent
For the basic theory of almost complex structures there are many references beyond the
seminal paper [NW] ([MDF] [MDF2]
[Si] in particular). Still, I shall often refer to [I-R] where special care has been taken in order
to present very simple proofs. [Ro2] may also be helpful
since it includes a long discussion of basic material
(using the notations to be introduced later, especially of $J$, $Q$ and the
Nijenhuis tensor). We refer to these sources for (standard) notions or notations that will not 
be explained below.
\bigskip\noindent
Pluripolarity in almost complex structures turns out to be a somewhat delicate matter.
Integrability near 0 in ${\bf R}^2$ of ${1\over r^2\log^2 r}$ ($r=\sqrt {x^2+y^2}$) leads
to our positive result, while non integrability of ${1\over r^2|\log r|}$ is crucial in the
counterexample.

% ppol0
\bigskip\bigskip

\noindent
{\bf 0. $\log$ and $\log |\log|$.}\hfill\break
In this section we work on ${\bf C}$, with variable $z$.
\medskip\noindent
One has ${\partial ^2\over \partial z\partial\overline z}
\log|z|=0$, on ${\bf C}\setminus \{ 0\}$.\hfill\break
Elementary computations give:
$${\partial^2 \over \partial z\partial\overline z}
\big( -\log|\log|z|\big)~={1\over 4 |z|^2\log^2|z|}~,$$
$${\partial^2 |z|\over \partial z\partial\overline z}
={1\over 4|z|}$$
(homogeneity being obvious).
\bigskip\noindent
As discussed in [Ro1] the function $-\log|\log|z||$ gives an
almost optimal spreading of mass for the Laplacian.
Almost complex structures on ${\bf C}^n$ are often treated, near a point
at which the structure is the standard one,
as deformations of the standard complex structure. This leads to study deformations
of the operator ${\partial^2 \over \partial z\partial\overline z}$
($={1\over 4} \Delta$). Corresponding to having the structure to be the
standard structure at the point, one is led to consider perturbations whose
terms of order 2 coincides with ${\partial^2 \over \partial z\partial\overline z}$
at 0.
\bigskip\noindent
Let $Q$ be a linear differential operator defined near 0, with say smooth coefficients
of the type
$$Q={\partial^2 \over \partial z\partial\overline z}+R_2(z,D)+Q_1(z,D)$$
with
$R_2(z,D)=\alpha_1(z){\partial^2\over \partial z\partial z}
+\alpha_2(z){\partial^2\over \partial z\partial\overline z}
+\alpha_3(z){\partial^2\over \partial \overline z\partial\overline z}$
and $\alpha_j(0)=0$,
$Q_1(z)=\beta_1(z){\partial\over\partial z}
       + \beta_2(z){\partial\over\partial \overline z}$, with no vanishing condition
of $\beta_j$ at 0.
\bigskip\noindent
From the estimates (for appropriate constant $C$)\hfill \break
$|\nabla (\log|\log|z||)|\leq {C\over |z||\log |z||}$ and
$|\nabla^2 (\log|\log|z||)|\leq {C\over |z|^2|\log |z||}$,
it is clear that near 0:
$$Q(-|\log|\log|z||)>0~.$$
That indicates a way that will be followed in section II.
\bigskip\noindent
There is no positivity of $\Delta (\log |z| )$ away from 0, so one
cannot hope to keep non negativity for $Q(\log |z|)$. A corrective term is
needed. Homogeneity considerations show that $C|z|^2$ cannot be enough,
but we have  $R_2(\log|z|)$ and $Q_1(\log |z|)=O({1\over |z|})$, 
$R_2(|z|)$ and $Q_1(|z|)=O(1)$. 
So, for $C>0$ large enough
$$Q (\log|z|~+C|z|)\geq 0~{\rm near}~0~.$$
That corresponds to the choice of the Chirka function, with logarithmic
pole at a point, on almost complex
manifolds, and it explains the need of having $J=J_{st}$ at the point.
\medskip\noindent
In case $R_2$ and $Q_1$ vanish to second, resp. first, order at 0, then for $C>0$ large enough
one can take the function $\log|z|+C|z|^2$ (or interestingly $\log|z|+C|z-z_0|^2$
for any fixed $z_0$), instead of $\log |z|+C|z|$. Part IV boils down to not much more
than that. 

% ppol1

\bigskip\bigskip

\noindent {\bf I. Preliminary Computations.}
Let $J$ be an almost complex structure defined on an open set $U\subset {\bf C}^n$.
We shall assume that $J$ is close enough to the standard complex structure $J_{st}$,
more precisely we assume that $J+J_{st}$ is invertible.\hfill\break
Then, as it is well known, 
the condition for $J$-holomorphy of a map $u$ from say the unit disc ${\bf D}$ 
in ${\bf C}$ (with coordinate $\zeta=x+iy$), into $U$ can be equivalently be given by:

\medskip
(a) the `Cauchy-Riemann' equation ${\partial u\over\partial y}=J(u){\partial u\over\partial x}$,
(where $J(u)$ is an ${\bf R}$-linear map satisfying $J^2=-{\bf 1}$).
\bigskip
\noindent or

(b) an equation 

$${\partial \overline u\over\partial \zeta}=Q(u){\partial u\over\partial \zeta},\qquad\quad (E1)$$
where $Q$ is a ${\bf C}$ linear operator (a complex $(n\times n)$ matrix valued function) of the
same smoothness as $J$, with $Q=0$ whenever $J=J_{st}$. The conjugate operator $\overline Q$
is given by $\overline Q=[J+J_{st}]^{-1}[J-J_{st}]$. A long discussion of that is in [Ro2].
\medskip\noindent
Let $u$ be a $J$-holomorphic disc, i.e. a map from ${\bf D}$ into $U$ satisfying (E1).
Let $\lambda$ be a real valued function defined on $U$, an elementary computation gives:
$${\partial^2\lambda\circ u\over\partial \zeta\partial\overline\zeta}~=~I+II+III\qquad\qquad (E2)$$
with
$$I=~\sum_{j,k\geq 1}^n{\partial^2\lambda \over\partial z_j\partial \overline z_k}
{\partial u_j\over\partial \zeta}{\partial\overline u_k\over\partial\overline\zeta}
+~\sum_{j,k\geq 1}^n{\partial^2\lambda \over\partial z_j\partial \overline z_k}
{\partial u_j\over\partial \overline \zeta}{\partial\overline u_k\over\partial\zeta}
$$
$$II=2~{\rm Re}~\sum_{j,k}{\partial^2\lambda\over\partial z_j\partial z_k}
{\partial u_j\over \partial\zeta} 
{\partial u_k\over\partial \overline \zeta}$$
$$III=2~ {\rm Re}~\sum_j {\partial \lambda\over\partial z_j}{\partial^2u_j\over\partial
\zeta\partial\overline\zeta}~~.$$
\bigskip\noindent
Standard holomorphic discs $u$ are harmonic (${\partial^2 u\over\partial\zeta\partial \overline \zeta}=0$),
$J$-holomorphic discs of course are not, but we still have an estimate of the Laplacian in terms of the
gradient. Indeed differentiation of (E1), with respect to $\overline \zeta$, gives:
$${\partial^2\overline
u\over\partial\zeta\partial\overline\zeta}
=
[Q_z(u).{\partial u\over \partial \overline\zeta}] {\partial u\over \partial \zeta}
+[Q_{\overline z}(u).{\partial \overline u\over \partial \overline\zeta}] {\partial u\over \partial \zeta}
+Q(u){\partial^2 u\over \partial \zeta\partial\overline\zeta}~.$$
Since
${\partial^2\overline u\over \partial \zeta\partial\overline\zeta}=
\overline{{\partial^2 u\over \partial \zeta\partial\overline\zeta}}$, if $Q(u)$ has an operator norm
$<{1\over 2}$, then one gets:
$$|{\partial^2u\over \partial \zeta\partial\overline\zeta}|
\leq 2~ |\big( [Q_z(u).{\partial u\over\partial \overline \zeta}]{\partial u\over\partial \zeta}
+[Q_{\overline z}(u).{\partial \overline u\over\partial \overline \zeta}]{\partial u\over\partial \zeta}
~\big)|.\qquad\qquad (E3)$$
Formula (E3) will be used, a first consequence that will not suffice for our purpose
is that, given bounds on $|\nabla Q|$, it gives an estimate $|{\partial^2 u\over \partial\zeta\partial\overline\zeta}|
\leq C|\nabla u|^2$.

% ppol2
\bigskip\bigskip

\noindent {\bf II. Pluripolarity of $J$-holomorphic curves.}
\medskip\noindent
After appropriate simple changes of variables (to straighten an embedded $J$-holomorphic disc 
and to make the almost complex structure standard along the straightened disc)
the heart of the matter is the following, on which we shall now concentrate
our efforts:
\bigskip\noindent
{\bf Proposition 1.} {\it In ${\bf C}^n$, with coordinates $Z=(z_1,Z')$ ($Z'=(z_2,\cdots ,z_n))$),
let $J$ be a ${\cal C}^2$ (or ${\cal C}^{1,1}$) almost complex structure
defined on a neighborhood of $\overline{\bf D}\times \{ 0\}$ (the set $|z_1|\leq 1$,
$z_2=\cdots =z_n=0$). Assume that $J=J_{st}$ along ${\bf C}\times \{ 0\}$. Then for $K>0$
large enough the function
$$\Lambda (Z)=-\log|\log |Z'||~+~K|z_1|^2$$
is $J$-plurisubharmonic on a neighborhood of $\overline {\bf D}\times \{ 0\}.$}
\bigskip\noindent
{\bf Proof.} 
Set $LL(Z')=-\log|\log |Z'||$.
The starting point is the estimate for the standard complex Hessian of the function
$LL(Z')$. See II.1 in [Ro1] 
$$4 \sum_{j,k\geq 2}{\partial ^2LL(Z')\over \partial z_j\partial \overline z_k}t_j\overline {t_k}
\geq {\sum_{j\geq 2}|t_j|^2\over |Z'|^2\log^2 |Z'|}~.$$
Let $u$ be a $J$-holomorphic disc (with values near ${\overline D}\times \{ 0\}$), we want to
show that \hfill\break
${\partial^2\over\partial\zeta\partial\overline \zeta}\Lambda\circ u~(0)\geq 0$.
\bigskip\noindent 
(1) We first estimate ${\partial^2
\over\partial\zeta\partial\overline \zeta}LL\circ u~ (0)$, by taking $\lambda
(Z)=LL(Z')$ in (E2). We shall estimate separately the terms I, II, and III, of course in the summations
only $j$ and $k\geq 2$ need to be considered. Various constants will be denoted by the
same letter $C$ in the proof.
\medskip\noindent
Set  $u(0)=(z_1^0,Z'_0)$,
$\tau = (\sum_{j=2}^n |{\partial u_j\over \partial \zeta}(0)|^2)^{1\over 2}$,
and $\epsilon = |{\partial u_1\over\partial \zeta}(0)|$. 
For the $\overline \zeta$ derivatives, for some appropriate
constant $C$,  one has the estimate
$|{\partial u_j\over \partial \overline \zeta}(0)|\leq C
|Z'_0|(\epsilon +\tau )$, according to (E1) and since $Q$ is of
class ${\cal C}^1$ and vanishes for $Z'=0$. It is an essential difficulty that
(E1) does not give a control of the $\overline \zeta$ derivatives of the $u_j$'s for $j\geq 2$
only in terms of the $\zeta$ derivatives of these functions alone, ${\partial u_1\over\partial \zeta}$
comes in the estimate. However, only the partial derivatives ${\partial u_j\over\partial \zeta}(0)$
for $j\geq 2$ contribute to the initial crucial positivity in the proof.
\bigskip\noindent
Estimate of I.
We have 
$$I\geq {\tau^2 \over 4 |Z'_0|^2\log^2(|Z'_0|)} ~.$$
Estimate of II. For an (other) appropriate constants $C$, $|\nabla^2LL(Z')|\leq 
{C\over |Z'|^2|\log |Z|'|}$.
So, for II we get (for some other appropriate constant $C$):
$$|II|\leq {C\over |Z'_0||\log |Z'_0||}( \epsilon +\tau )\tau~.$$
Estimate of III. One first has an estimate $|\nabla  LL(Z')|\leq {C\over |Z'||\log |Z'||}.$
More delicate is the estimate of ${\partial^2u\over \partial\zeta\partial\overline \zeta}(0)$.
Here we need to analyze both terms on the right hand side of (E3).
Using that $Q=0$ if $Z'=0$, that gives $|{\partial u\over\partial\overline\zeta}(0)|\leq C |Z'_0|
|{\partial u \over \partial \zeta}(0)|=C|Z'_0| (\epsilon +\tau )$, one gets
$$[Q_z(u).{\partial u\over\partial \overline \zeta}]{\partial u\over\partial \zeta}(0)\leq
C|Z'_0|(\epsilon^2 +\tau^2)~.$$
Next we evaluate (for $\zeta=0$) $Q_{\overline z}(u).\overline{\partial u\over\partial \zeta}=
a+b$, where $a=\overline{\partial u_1\over\partial \zeta} {\partial Q\over\partial \overline z_1}$
and $b=\sum_{j=2}^n \overline{\partial u_j\over\partial \zeta}{\partial Q\over\partial\overline z_j}$.
Since $Q\equiv 0$ for $Z'=0$, one has $|{\partial Q\over\partial \overline z_1}(Z)|\leq C |Z'|$, and this is
where we need ${\cal C}^2$ regularity of $Q$ (or at least ${\cal C}^{1,1}$ regularity).
Then one gets $|a|\leq C |Z'_0|\epsilon $. For the $b$ term, one sees that
$|b|\leq C\tau $. Therefore
$|[Q_{\overline z}(u).\overline{\partial u\over\partial \zeta}]{\partial u\over\partial\zeta}(0)|
\leq C~(|Z'_0|\epsilon +\tau) (\epsilon +\tau))$.
Consequently
$$|{\partial^2u\over \partial \zeta\partial\overline \zeta}(0)|\leq
C |Z'_0|(\epsilon^2 +\tau^2)~+~C\tau (\epsilon +\tau )~.$$
Finally
$$|III|\leq C ({\epsilon^2 \over |\log |Z'_0||}+{\tau^2\over |Z'_0||\log|Z'_0||}+
{\epsilon\tau \over |Z'_0||\log|Z'_0||})~.$$
\bigskip\noindent
Putting together all the estimates obtained so far,
we get (for new constants $C$):
$${\partial^2\over\partial\zeta\partial\overline\zeta}LL\circ u~(0)\geq
{\tau^2\over 4|Z'_0|^2\log^2|Z'_0|}$$
$$-{C\tau^2\over |Z'_0||\log|Z'_0||}-{C\epsilon\tau\over  |Z'_0||\log|Z'_0||}
-{C\epsilon^2\over |\log|Z'_0|}~.$$
\smallskip\noindent
On the right hand side, the first term gives positivity, and the worst term is
$-{C\epsilon \tau\over |Z'_0||\log |Z'_0||}$.
\bigskip\noindent
Remark: Before finishing the proof it is good to point out that if $\epsilon =0$,
the right hand side
is clearly non negative provided $|Z'_0|$ is small enough, i.e. $u(0)$ is close
to ${\bf C}\times \{ 0\}$. In fact the estimates above prove that in ${\bf C}^N$ equipped with
an almost structure $J$ of class ${\cal C}^1$, with $J(0)=J_{st}$, the function
$-\log|\log |Z|$ is $J$-plurisubharmonic near 0. No additional corrective term is needed.
The proof is just a matter of dropping
all the terms with index $j=1$ above. In particular the term $a$ in the estimate of III above
need not be considered and this is why ${\cal C}^1$ smoothness would be enough.
Also, note that at any rate $J$-purisubharmonicity of
$-\log |\log |Z||$ is, for ${\cal C}^2$ data, a completely trivial
consequence of the Proposition, by consideration of ${\bf C}^N\times {\bf C}_{st}$.
\bigskip
\noindent
(2) We now estimate ${\partial^2\over\partial\zeta\partial\overline \zeta}|u_1|^2$.
So this time we apply (E2) to the function $\lambda (Z)= |Z_1|^2$.
With the notation as in (1), one gets immediately
$$I\geq \epsilon^2~~~~,~~II=0~.$$
For III, ${\partial \lambda\over \partial z_1}=\overline z_1$, 
and  ${\partial \lambda\over \partial z_j}=0$ if $j\neq 1$, but we need to estimate
${\partial^2u_1\over \partial\zeta\partial\overline\zeta}$ as
in (1). This gives us
$$|III|\leq C\big( |Z'_0|\epsilon^2+\tau^2+\epsilon \tau)~~.$$
\bigskip\noindent
(3) According to (1) and (2) $\Lambda=LL(Z')+K|Z_1|^2$ will be $J$-plurisubharmonic 
near $\overline{\bf D}\times \{ 0\}$ if  for $|Z'_0|$ small enough, and for any $\epsilon$ 
and $\tau\geq 0$,
$$
{\tau^2\over 4|Z'_0|^2\log^2|Z'_0|}+ K\epsilon^2$$
$$
-{C\tau^2\over |Z'_0||\log|Z'_0||}-{C\epsilon\tau\over  |Z'_0||\log|Z'_0||}
-{C\epsilon^2\over |\log|Z'_0||}-CK|Z'_0|\epsilon^2-CK\tau^2-CK\epsilon\tau\geq 0~.$$
By taking $K>4 C^2$, one has
$${\tau^2\over 4|Z'_0|^2\log^2|Z'_0|}+ K\epsilon^2-{C\epsilon\tau\over  |Z'_0||\log|Z'_0||}
\geq {1\over 8} \big( {\tau^2\over |Z'_0|^2\log^2|Z'_0|}+ K\epsilon^2\big)~.$$
Once $K$ has been fixed, the other terms are immediately absorbed when
$|Z'_0|$ is small enough. This ends the proof of the Proposition.

\bigskip\noindent
Remark: To get a logarithmic singularity, one would be tempted to add 
$c\log |Z'|$ to $\Lambda$, for $c>0$ small enough. Then, coming from term II in E2, and from estimates
of derivatives,
a term arises
with ${1\over |Z'_0|}$ instead of ${1\over |Z'_0||\log |Z'_0||}$. As shown in Part III,
this term cannot be controlled.

% ppol 3

\bigskip\bigskip

\noindent{\bf III. An Example.}
\hfill\break
On ${\bf C}^2$, we consider the almost complex structure $J$ defined by
taking
$$Q=\left(\matrix {0 & 0 \cr
                 z_2 & 0 \cr}\right)~~.$$
So the equations for $J$-holomorphic discs $\zeta\mapsto u(\zeta )=\big( u_1(\zeta ),u_2(\zeta )\big)$
are: $u_1$ must be holomorphic, and  
${\partial \overline u_2\over\partial \zeta}=u_2{\partial u_1\over\partial \zeta}$.
\bigskip\noindent
It then follows that ${\partial u_2\over\partial \overline \zeta}=
\overline u_2\overline {\partial u_1\over\partial \zeta}$, and
${\partial^2u_2\over\partial\zeta\partial\overline\zeta}=u_2
|{\partial u_1\over\partial\zeta}|^2$ .
\medskip\noindent
The fact that $Q_z\neq 0$ is important, this is related to the non vanishing
of the Nijenhuis tensor. See detailed discussion in [Ro2] (Proposition 6), and IV below.
A direct description of $J$ is given at the end of Part III.
\medskip\noindent
Note that $J=J_{st}$ along ${\bf C}\times \{ 0 \}$, that is therefore
$J$-holomorphic.
\bigskip\noindent
{\bf Proposition 2.} {\it There is no $J$-plurisubharmonic function $\lambda$,
defined near 0 and not identically $-\infty$ near 0, such that
$\lambda (z_1,z_2)\leq \log |z_2|$.}
\bigskip\noindent
Recall that, for plurisubharmonic functions, we restrict our attention to
functions that are ${\cal C}^2$ off their $-\infty$ set. I have not tried
to extend Proposition 2 to the case of merely continuous plurisubharmonic
functions. The proof will be by contradiction, so we assume that there
is a function $\lambda$ with the desired properties.
\bigskip\noindent
We rewrite equation (E3), that simplifies, and with a different arrangement of terms,
Let $u=(u_1,u_2)$ be a $J$-holomorphic map from a neighborhood of 0 in ${\bf C}$ into some
fixed neighborhood of 0 in ${\bf C}^2,J)$. \medskip
$${\rm Set}~u(0)=(z_1,z_2)~~.~{\rm Assume}~{\partial u_1\over \partial \zeta} (0)=1~~,
~{\rm set}~ {\partial u_2\over \partial \zeta} (0)=t~.$$ 
\smallskip\noindent
Note that except if ${\partial u_1\over \partial \zeta} (0)=0$,
${\partial u_1\over \partial \zeta} (0)=1$ can be achieved by simple rescaling.
By the equations for $J$-holomorphy: ${\partial u_1\over\partial\overline \zeta}=0$,
and ${\partial u_2\over\partial\overline \zeta}(0)= \overline z_2$.
One gets
$${\partial^2\over\partial\zeta\partial\overline\zeta}(\lambda\circ u)(0)=A+B+C_1+C_2 \qquad\qquad (E4)$$
with:
$$A= 2~{\rm Re}[({\partial^2\lambda \over \partial \overline z_1\partial z_2}
+ {\partial^2\lambda \over \partial z_2^2}\overline z_2)t]~,$$
$$B= {\partial^2\lambda\over \partial z_2\partial\overline z_2}~(|t|^2+|z_2|^2)~,$$
$$C_1= 2~{\rm Re}[({\partial^2\lambda \over \partial z_1\partial z_2} \overline z_2)
+ ({\partial\lambda \over \partial z_2} z_2)]~,$$
$$C_2= {\partial^2\lambda\over \partial z_1\partial\overline z_1}~
~.$$
\bigskip\noindent
The strategy for proving the Proposition will be to find $K_1>0$ 
such that for $z_1=0$ and some $z_2$'s, with $|z_2|$ arbitrarily small, one can choose
$t$  such that $A+B+C_1+C_2\leq K_1\log|z_2|~<0$. Negativity will come from
the $A$ term, more precisely from ${\partial^2\lambda \over \partial z_2^2}\overline z_2$, coming from II
which was previously treated as an error term. In our example, one can avoid cancellation by
${\partial^2\lambda\over \partial\overline z_1\partial z_2}$  by appropriate choice of the argument of $z_2$.
The other terms will be `error terms'.
\hfill\break It happens that a good choice of $t$ will be with
$|t|=K_2|z_2|~|\log |z_2||$, for some appropriate value of $K_2$.

\bigskip\noindent
{\bf Preliminary Remarks.}\hfill\break
For any fixed $z_1\in {\bf C}$, $\zeta\mapsto (z_1,\zeta )$ is a $J$-holomorphic
map. Hence $z_2\mapsto \lambda (z_1,z_2)$ is a subharmonic function.
Because of the Logarithmic singularity we must have
$$\lambda (z_1,z_2)=a(z_1)\log |z_2|~+~\mu (z_1,z_2)~,\qquad\qquad (E4)$$
where $a(z_1)\geq 1$ (the Lelong number) and the measure $\nu_{z_1}$,
on ${\bf C}$ (near 0), defined by 
$\nu_{z_1}~={\partial^2\over\partial z_2\partial \overline z_2}
(\mu (z_1,.)~)$ has no point mass at $z_2=0$. So, $\mu (z_1,z_2)$ can be thought
as a term of lesser singularity along $z_2=0$.
\bigskip\noindent
{\bf 1.} We shall first prove the non existence of $\lambda$ if one imposes some
conditions that seem very natural and reasonable, in view of the above discussion.
\bigskip\noindent
The hypotheses that we now make on the term $\mu$,`less singular' than
$\log|z_2|$, are:
\bigskip\noindent
For $z_1$ small, and as $z_2\to 0$:
$$z_2^2{\partial^2\mu\over \partial z_2^2}\to 0~~,~~z_2{\partial\mu\over \partial z_2}\to 0~~.
\qquad \qquad ({\bf H})$$
Since $\mu$ is real valued the second condition implies that ${\mu\over\log|z_2|}\to 0$,
but we won't need this remark.
\medskip\noindent
Since we clearly have invariance under $z_1$ translations, we can partially smooth $\lambda$
(and thus $a$ and $\mu$) by convolution in the $z_1$ variable. Then $a(.)$ is smooth and for fixed $z_2$, a sup norm
estimate of $\lambda$ or $\mu$  or of some derivative of these functions, results in a similar sup norm
estimate for functions obtained by further differentiations
in the $z_1$ direction. We can therefore assume that additionally
to $({\bf H})$, we also have, for some constant $M>0$:
$$|{\partial^2\lambda \over \partial z_1\partial \overline z_1}|\leq M~|\log|z_2||~~,~~
|{\partial^2\lambda \over \partial z_1\partial z_2}|\leq {M\over |z_2|}~~,~~
z_2{\partial^2\mu\over \partial \overline z_1\partial z_2}\to 0,~{\rm as}~z_2\to 0~.\quad ({\bf H^+})$$
\bigskip\noindent
Proof of Proposition assuming $({\bf H})$ (and thus $({\bf H^+})$.
\medskip\noindent
We first evaluate $A+C_1+C_2$. We have
$${\partial^2\lambda\over \partial \overline z_1\partial z_2}=
{1\over 2} {\partial a\over\partial \overline z_1}{1\over z_2}
+{\partial^2\mu \over \partial \overline z_1\partial z_2}~=~
{1\over 2} {\partial a\over\partial \overline z_1}{1\over z_2}~+~o({1\over |z_2|})~,$$
$${\partial^2\lambda \over \partial z_2^2}={-a(z_1)\over 2 z_2^2}~+~o({1\over |z_2|^2})~~.$$
Consequently
$${\partial^2\lambda\over \partial \overline z_1\partial z_2}
+ {\partial^2\lambda \over \partial z_2^2}\overline z_2=
{1\over 2}\big( {\partial a\over\partial \overline z_1}-a(z_1){\overline z_2\over z_2}\big)
{1\over z_2}~+~o({1\over |z_2|})~.$$
For any fixed $z_1$, and we take $z_1=0$, there exits an interval ${\cal I}$ in
${\bf R} / 2\pi {\bf Z}$,
of length $|{\cal I}|\geq {\pi\over 4}$, such that if $z_2=re^{i\theta}$ with $\theta\in {\cal I}$,
(avoiding cancellation of $-a$ by ${\partial a\over \partial\overline z_1}$) on has
 $$|{\partial^2\lambda\over \partial \overline z_1\partial z_2}
+ {\partial^2\lambda \over \partial z_2^2}\overline z_2|\geq {1\over 3 |z_2|}~~,$$
if $|z_2|$ is small enough.
\medskip\noindent
Under our hypotheses, as $z_2$ approaches 0, the $C_1$ term stays bounded and for the 
$C_2$ term one gets a bound $|C_2|\leq M|\log|z_2||$.
Therefore if we take $K$ large enough, for any $z_2=|z_2|e^{i\theta}$, with $\theta\in {\cal I}$,
we can choose $t$ with $|t|=K|z_2||\log|z_2||$, with the argument chosen so that
$$A+C_1+C_2\leq -{K\over 4}~|\log|z_2||~.$$
\bigskip\noindent
We finally need to control the $B$ term. Ideally we would like 
$(|z_2|^2\log|z_2|){\partial^2\lambda\over \partial z_2\partial\overline z_2}\to 0$.
That would make the $B$ term negligible. Note that this is asking more than the natural 
requirement suggested by `homogeneity' consideration, that would be
$(|z_2|^2){\partial^2\lambda\over \partial z_2\partial\overline z_2}\to 0$,
since ${\partial^2\lambda\over \partial z_2\partial\overline z_2}={
\partial^2\mu\over \partial z_2\partial\overline z_2}$.
\medskip\noindent
Since, as pointed out earlier, $\lambda$ is a subharmonic function of
$z_2$ (for $z_1$ fixed), its Laplacian is a bounded positive measure.
For fixed $z_1$, ${\partial^2\lambda \over \partial z_2\partial \overline z_2}$
is (near $z_2=0$) the sum of a point mass at 0 and of an integrable function.
So, given any $\epsilon >0$ the inequality 
$|{\partial \lambda\over \partial z_2\partial\overline z_2}| > {\epsilon\over
|z_2|^2|\log|z_2||}$ cannot be satisfied near 0 by all $z_2$ in a given sector.
Then, at any point $z_2$ where the reverse inequality 
$|{\partial \lambda\over \partial z_2\partial\overline z_2}|\leq {\epsilon\over
|z_2|^2|\log|z_2||}$ 
holds, for $t$ as above one has
$$|B|=|{\partial^2\lambda\over\partial z_2\partial\overline z_2}(|t|^2+|z_2|^2)|
\leq {\epsilon\over |z_2|^2|\log|z_2||}K^2|z_2|^2\log^2|z_2| + {\epsilon\over |\log|z_2||}
\leq\epsilon K^2|\log|z_2||~+1.$$
By taking $\epsilon$ so that $\epsilon K^2<{K\over 8}$, one gets
$A+B+C_1+C_2\leq - {K\over 8}|\log |z_2||~+1<0$ ,($|z_2|$ small). Therefore $\lambda$ is not $J$-plurisubharmonic since its 
restriction to any $J$-holomorphic disc satisfying $u(0)=(0,z_2)$, ${\partial u_1\over \partial \zeta}(0)=1$
and ${\partial u_2\over \partial \zeta}(0)=t$, will not be subharmonic.

% ppol32

\bigskip\bigskip

\noindent {\bf 2.} We now prove Proposition 3 in full generality.
Our first task is to get that $z_2{\partial\mu\over\partial z_2}$
tends to 0, uniformly with respect to $z_1$, as $z_2$ approaches zero
along some `fat' set.
\bigskip
\noindent
A subset $E\subset {\bf C}$ will be said to be {\it fat} if 0 is a point of Lebesgue density
of $E$, i.e. 

$${|E\cap\{ |z|<r\}|\over \pi r^2}\to 1~~~{\rm as}~r\to 0~.$$
\bigskip\noindent
Assume that there is a function $\lambda_1$ with the properties mentioned in the statement
of the Proposition, defined for $|z_1|\leq 4R$, $|z_2|\leq 3R$. For $0< |z_2|<2R$ (so $z_2\neq 0$),
set $$\nu (z_2)=\int_{|z_1|<3R} {\partial^2 \lambda_1\over \partial z_2\partial\overline z_2}(z_1,z_2)~dm(z_1)$$
where $dm$ denotes Lebesgue measure. By subharmonicity of $\lambda_1$ in the variable $z_2$,
$\nu$ is an integrable function. Let $\lambda$ be a function
obtained by partially smoothing $\lambda_1$ by convolution in the $z_1$ variable,
with a non non negative function of $z_1$, of small support.
$J$-plurisubharmonicity is preserved since $J$ is invariant under $z_1$ translations,
and we replace the original function $\lambda_1$ by $\lambda$.
Then there is $C>0$ such that for all $|z_1|<2R$ we have
$$0\leq {\partial^2 \lambda\over \partial z_2\partial\overline z_2}(z_1,z_2)\leq C \nu (z_2)~.$$
Let $\mu$ be the function defined in (E4). 
Let $\chi (z_2)$ be the characteristic function of the set $\{ |z_2|< 2R\}$.
Then for each fixed $z_1$, on the set $|z_2|<R$,
$\mu$ is the sum of a smooth function
(with uniform bounds) and of 
$(\chi {\partial^2\mu\over\partial z_2\partial\overline z_2})\ast 
{2\log|z_2|\over \pi}$  
(convolution, in the $z_2$ variable, of the Laplacian with the Newtonian potential). 
And
${\partial ^2\mu \over \partial z_2\partial \overline z_2}
$ is simply the restriction to $z_2\neq 0$ of 
${\partial \lambda\over\partial z_2\partial\overline z_2}$ (the point mass at $\{ 0\}$ being
thus dropped). So, modulo smooth function
${\partial \mu\over \partial z_2}=
(\chi {\partial^2\mu \over \partial z_2\partial\overline z_2})\ast {1\over \pi z_2}$. 
We therefore get that for every $|z_1|<2R$, for $|z|_2|<R$

$$|{\partial \mu\over\partial z_2}|\leq C{1\over |z_2|}\ast \nu~+C$$
where the convolution takes place in the $z_2$ variable.
It then follows from Lemma A1 in the appendix that there is a fat set $E$ such that:
$${\rm Sup}_{|z_1|\leq R} |z_2{\partial \mu\over\partial z_2}(z_1,z_2)|~\to 0~,~{\rm as}~z_2\to 0~,z_2\in E~.$$
So we have a replacement of the second hypothesis in ${\bf H}$, if we restrict $z_2$ to $E$.
\bigskip
\noindent
It is essential that for $z_2\in E$ we had estimates, uniformly for  
$|z_1|<R$. As done in ${\bf 1}$ we can replace $\lambda$ by a function obtained
by a second (probably un-needed) smoothing by convolution in the $z_1$ variable and
if one restricts to $z_2\in E$, hypotheses ${\bf H^+}$  are satisfied.
\bigskip\noindent
We finally need a replacement for the first hypothesis in ${\bf H}$. This is a more subtle matter
since one has to consider now a singular integral. 
$${\partial^2\mu\over\partial z_2^2}~=~\nu\ast {-1\over \pi z_2^2}~+g,$$
where $g$ is a function bounded near 0.
\bigskip\noindent
For the end of the proof, contrary to above, we have to fix $z_1$. Take $z_1=0$. Lemma A2 in the Appendix
shows that there exists a fat set $E'$ ($\subset E$) on which additionally to the above we have
$z_2^2{\partial^2\mu\over\partial z_2^2}\to 0$ as $z_2\to 0$, $z_2\in E'$.
\bigskip\noindent
Lemma A3 asserts that for any $\epsilon >0$ the inequality ${\partial^2\lambda\over\partial z_2\partial
\overline z_2}> {\epsilon\over |z|^2\log|z|}$ cannot hold for all $z_2$ near 0 $z_2\neq 0,~{\rm and}~
z_2\in E'$. From here the proof is identical to the proof in ${\bf 1}$, with the only difference that
one has to take $z_2\in E'$.
\bigskip\bigskip

% ppol33

\noindent
{\bf 3.} Tensor $J$. The almost complex structure of the example
on ${\bf C}^2$ has been given by the ${\bf C}$-linear operator $Q$ that characterizes
$J$-holomorphy:
$${\partial\overline u\over\partial\zeta}=Q(u){\partial u\over\partial\zeta}~.$$
We now wish to write the corresponding ${\bf R}$ linear operator $J$ 
defining the almost complex structure on ${\bf R}^4$,
on which coordinates will be $(x_1,y_1,x_2,y_2)$ (identified with ${\bf C}^2$,
by setting $z_j=x_j+y_j$).
\medskip\noindent
There is a formula:
$J=J_{st}[{\bf 1}-\overline Q][{\bf 1}+\overline Q]^{-1}$ (see 2.1 in [Ro2]) that
one could use. We do a direct computation instead.
We observe that $\zeta =(x+iy)~\mapsto (z_1+\zeta,z_2+\overline z_2\overline\zeta )$ satisfies the equation for
$J$-holomorphicity at $\zeta =0$. In terms of $J$ the equation is
${\partial u\over \partial y}=J(u){\partial u\over \partial x}$.
Therefore
$$[J(z_1,z_2)] (1,\overline z_2)=(i,-i\overline z_2)~.$$
Since $J^2=-{\bf 1}$, $[J(z_1,z_2)] (i,-i \overline z_2)=(-1,-\overline z_2)~.$ \hfill\break
Since $\zeta\mapsto (z_1,z_2+\zeta)$ is $J$-holomorphic, one also has
$$[J(z_1,z_2)] (0,\xi )=(0,i\xi )~, {\rm for~any}~\xi\in {\bf C}~.$$
One immediately gets  
$$[J(z_1,z_2)](1,0)=(i,-2i\overline z_2)~,~
[J(z_1,z_2)](i,0)=(-1,-2\overline z_2)~,$$
$$[J(z_1,z_2)](0,1)=(0,i)~,
~[J(z_1,z_2)](0,i)=(-1,0)~.$$
With real ${\bf R}^4$ notations at point $(x_1,y_1,x_2,y_2)$, $J$ is given
by the matrix
$$ J=    \left(\matrix{ 0  &   -1  &  0 & 0      \cr
                         1  &   0  &   0 & 0     \cr
                         -2y_2 & -2x_2 & 0 & -1  \cr
                         -2x_2  &   2y_2 & 1 & 0     \cr}\right) \quad .$$

% ppol4
\bigskip\bigskip

\noindent
{\bf IV. Logarithmic singularities.}
\hfill\break
There is no surprise that if we assume better matching with the standard complex
structure along an embedded $J$-holomorphic curve, there will be a $J$-plurisubharmonic
function with logarithmic singularity along that curve. Better matching with the standard
complex structure is given by vanishing of the Nijenhuis tensor.
\bigskip\noindent
On an almost complex manifold $M$ a $(0,1)$ vector field is a complexified tangent vector field 
$\overline L$ such that $\overline L =X+iJX$ for some real tangent vector field.
The vanishing of the Nijenhuis tensor at a point $p\in M$ is equivalent to the fact
that if $\overline L_1$ and $\overline L_2$ are $(0,1)$ vector fields, their Lie bracket
is of type $(0,1)$ at $p$ (i.e. $[\overline L_1,\overline L_2](p) = X_p+iJ(p)X_p$, for some real tangent vector 
$X_p$ to $M$ at $p$). Although smoothness hypotheses are made precise below, we did not try to improve them,
e.g. by using the Whitney extension Theorem.
\bigskip\noindent
{\bf Lemma.} {\it Let $J$ be an almost complex structure of class ${\cal C}^4$
defined on a neighborhood
of $\overline {\bf D}\times \{ 0\}$ in ${\bf C}^n$. Assume that $J=J_{st}$
along ${\bf C}\times \{ 0\}$. The following are equivalent:
\medskip
(a) The Nijenhuis tensor $N_J$ vanishes along ${\bf D}\times \{ 0\}$
\medskip
(b) There exists a ${\cal C}^3$ local change of coordinates given by:
$$Z_r(z)=z_r+\sum_{k,l\geq 2}a^r_{k,l}(z_1)z_k\overline z_l +
        +\sum_{k,l\geq 2}b^r_{k,l}(z_1)\overline z_k\overline z_l~,$$
such that for any ${\cal C}^2$ $(0,1)$ vector field $\overline L$
$$\overline L(Z_j)(z_1,z_2,\cdots ,z_n)~=
~O(\sum_{j\geq 2}|z_j|^2)~,$$
in a neighborhood of ${\bf D}\times \{ 0\}$.}
\bigskip\noindent 
The proof is identical to the proof of Proposition 4 in [Ro2], with $z_1$
as a parameter. In the new coordinates the almost complex
structure is of class ${\cal C}^2$ at least.
We shall prove that (a) implies (b). The converse is 
easy since $\nabla Q$, hence $\nabla J =0$ in the new coordinates, and evaluating
the Nijenhuis tensor requires only one derivative of $J$. \hfill\break
{\bf Proof:} By simple linear algebra there is a basis of $(0,1)$ vector fields
$\overline L_j$ ($j=1,\cdots ,n)$ such that
$$\overline L_j={\partial\over\partial\overline z_j}+
\sum_q\alpha_{j,q}(z){\partial\over\partial z_q}~,$$
where the functions $\alpha_{j,q}$ vanish along ${\bf D}\times \{ 0\}.$
(a) reduces to the simple form $[\overline L_j,\overline L_k]=0$
along ${\bf D}\times \{ 0\}$, giving us
$${\partial \alpha_{j,q}\over \partial \overline z_k}(z_1,0,\cdots )=
{\partial \alpha_{k,q}\over \partial \overline z_j}(z_1,0,\cdots )~.$$
\hfill\break
We have
$$\overline L_j(Z_r)=\alpha_{j,r}(z)+\sum_{k\geq 2}a^r_{k,j}(z_1)z_k
+\sum_{k\geq 2}(b^r_{k,j}(z_1)+b^r_{j,k}(z_1) )\overline z_k~+O(\sum_{j\geq 2}|z_j|^2)~.$$
We then choose $a^r_{k,j}(z_1)=-{\partial \alpha_{j,r}\over\partial z_k}(z_1,0, \cdots )$
and $b^r_{k,j}(z_1)$ such that
$b^r_{k,j}(z_1)+b^r_{j,k}(z_1)=-{\partial \alpha_{j,r}\over\partial\overline z_k
}(z_1,0\cdots )$. This is possible since 
${\partial \alpha_{j,q}\over \partial \overline z_k}(z_1,0,\cdots )=
{\partial \alpha_{k,q}\over \partial \overline z_j}(z_1,0,\cdots )~.$
\bigskip\noindent
In the coordinates $(Z_j)$ of the Lemma, $J$ coincides to $J_{st}$ to second order
along ${\bf D}\times \{ 0\}$. So, along ${\bf D}\times \{ 0\}$, not only
$Q=0$ but also $\nabla Q=0$. Recall that $Q_{\overline z}=0$
can always be achieved according to Sukhov and Tumanov (see [D-S] Lemma 3.2, sections 2.2 in [S-T1 ] [S-T2 ],
and [Ro2] in particular Proposition 6).
It is for getting $Q_z=0$ that one needs the vanishing of the Nijenhuis tensor.
\bigskip\noindent
{\bf Proposition 3.} {\it In ${\bf C}^n$, with coordinates $Z=(z_1,Z')$ ($Z'=(z_2,\cdots ,z_n)$),
let $J$ be a ${\cal C}^4$ almost complex structure
defined on a neighborhood of $\overline{\bf D}\times \{ 0\}$ (the set $|z_1|\leq 1$,
$z_2=\cdots =z_n=0$). Assume that 
along ${\bf C}\times \{ 0\}$, $J=J_{st}$ and the
Nijenhuis tensor $N_J$ vanishes. Then for $K>0$
large enough the function
$$\Lambda (Z)=\log |Z'|~+~K|Z|^2$$
is $J$-plurisubharmonic on a neighborhood of $\overline {\bf D}\times \{ 0\}.$}
\bigskip\noindent
{\bf Proof.} The proof is similar to the proof in part {\bf II}, and simpler.
Due to the second order vanishing of $Q$ along ${\bf C}\times \{ 0\}$.
Quickly said:
$\log |Z'|\circ u$ (that would be subharmonic if $u$ were holomorphic)
has a Laplacian bounded from below, if $u$ is $J$-holomorphic, and $|\nabla u|\leq 1$. Then
$|Z|^2$ brings enough positivity.\hfill\break
We adopt the same notations as in the Proof of Proposition 3.
We estimate 
${\partial^2(\log|z'|\circ u)\over\partial\zeta\partial\overline \zeta}(0)$ by taking
$\lambda =\log |Z'|$ in $(E2)$.
We have $I\geq 0$ since $\log |Z|$ is $J_{st}$-plurisubharmonic.
\smallskip\noindent
For II $|{\partial u_k\over
\partial \overline \zeta}|(0)\leq C |Z'_0|^2|(\epsilon +\tau )$ (with now $|Z'_0|^2$
instead of $|Z'_0|$), but now $|\nabla^2\lambda|\leq {C\over |Z'_0|^2}$, so
(with a new constant $C$) $II\leq C(\epsilon +\tau )\tau$. 
\smallskip\noindent
For III, there is a gain in estimating the $b$ term in the proof of Part {\bf II},
we now have ${\partial Q\over \partial \overline z_j}=0$ for $Z'=0$, not only for
$j=1$. This results in having
$|{\partial^2u\over \partial \zeta\partial\overline \zeta}(0)|\leq
C |Z'_0|(\epsilon^2 +\tau^2)~.$ Thus $|III|\leq C (\epsilon^2 +\tau^2)$
and finally
$${\partial^2\log (|Z'|\circ u)\over
\partial \zeta\partial\overline\zeta}(0)\geq -C(\epsilon^2 +\tau^2 )~.$$
\smallskip\noindent
If $Z'_0$ is close to 0, one has 
${\partial^2(|Z|^2\circ u)\over \partial\zeta\partial\overline \zeta}(0)\geq M (\epsilon^2+\tau^2)$
($M$ a positive constant). This can be obtained by the same considerations as above 
or more simply by using that
along ${\bf C}\times \{ 0\}$ $dd^c_J=dd^c$ and the simple formula 

$$\Delta (\lambda\circ u)(\zeta)=
[dd^c_j\lambda ]_{u(\zeta )} 
\big({\partial u\over\partial x},[J(u)]{\partial u\over\partial x})\big) ,~$$ see e.g. Lemma 1.2 in [I-R].
\smallskip\noindent
If $K$ is large enough, one therefore gets ${\partial^2 (\Lambda\circ u)
\over\partial \zeta\partial\overline\zeta}(0)\geq 0$. Proposition 3
is thus proved.
 
% ppol5
\bigskip\bigskip

\noindent{\bf V. Non-pluripolarity of generic submanifolds.}
\bigskip\noindent
It is often the case that an argument using Jensen measures can be avoided
and can be replaced by an elementary reasoning. However (at least for people
with some training in function algebras) Jensen measures give an immediate
insight. So we start by discussing Jensen measures although we shall say later
how elementary reasoning suffices (so {\bf 1} can be skipped).
\bigskip\noindent
{\bf 1) Jensen Measures.} No originality can be claimed for the following:
\medskip\noindent
{\bf Lemma 2.}{\it Let $K$ be a compact set in an almost complex manifold 
$M$, and let $p\in M$. Assume that for every $J$-plurisubharmonic function $\lambda$,
$\lambda (p)\leq {\rm Sup}_K\lambda$. Then there exists a probability measure
$\mu$ on $K$ such that for very $J$-plurisubharmonic function $\lambda$,
$\lambda (p)\leq \int_K\lambda~d\mu$.}
\medskip\noindent
Such a measure will be called a Jensen measure for $p$ on $K$.
The proof is the standard proof of the existence of Jensen measures
in the theory of uniform algebras (see e.g. [Ga],or [St] page 11).

\medskip\noindent
{\bf Proof.} 
Since from the beginning we have restricted our attention to $J$-plurisubharmonic
functions $\lambda$ that are continuous off their $-\infty$ set, we can restrict our attention
to those that are continuous and that take values in $(-\infty ,+\infty)$,
approximating $\lambda$ by $\chi_N\circ\lambda$, where $\chi_N$ is a convex increasing function
defined on ${\bf R}$ that satisfies $\chi_N(t)=-N$ if $t\leq -N-1$ and
$\chi_N(t)=t$ if $t\geq -N+1$,
with $N\to +\infty$. 

\smallskip\noindent
Let $\Gamma$ be the open convex cone in the space of real valued continuous functions
on $K$ that consists of the functions $\varphi$ for which there is a $J$-plurisubharmonic
function $\lambda$ defined on $M$, such that $\lambda (p) >0$ and $\lambda \leq \varphi$ on $K$.
Then $0\notin \Gamma$. By the Hahn-Banach Theorem, there exists a linear functional $\phi$
on the space of continuous functions on $K$, that
is positive on $\Gamma$. All positive functions on $K$ belong to $\Gamma$ (by taking
constant functions for $\lambda$), so $\phi$ is a positive functional.
After multiplication by appropriate constant, $\phi$ is given by a probability measure $\mu$.
If $\alpha \in {\bf R}$ and $\lambda$ is a $J$-plurisubharmonic function such that
$\alpha <\lambda (p)$, then $(\lambda - \alpha) \in \Gamma$, so
$\int_K(\lambda -\alpha )~d\mu >0$. Hence $\alpha < \int_K\lambda$ since $\int_Kd\mu =1$.
Q.E.D. 
\bigskip\noindent
{\bf 2) Dimension 2.}
First, we show that if $\lambda$ is a (standard) plurisubharmonic function defined near 0 in ${\bf C}^2$
that is $-\infty$ on ${\bf R}^2$ near 0, then $\lambda$ is identically $-\infty$ on ${\bf C}^2$
near 0. There are of course many immediate arguments for showing that. The one here is not the simplest
one (by far) but the advantage if that it extends immediately to the almost complex setting,
by treating the almost complex case, just as a small perturbation of the complex case, 
as it is usual to do, with no delicate analysis needed this time.
\bigskip\noindent
For $r>0$ small enough let $U=\{
(z_1,z_2),~ {\rm Im}~z_1>0,~|{\rm Im}~ z_2|\leq ({\rm Im}~z_1)^2~,~|z_1|^2+|z_2|^2<r^2\}$.
The boundary of $U$ consists of 3 parts $A,B$ and $C$: with $A\subset {\bf R}^2$,
$B$ a subset of the sphere of radius $r$, and $C\subset \{ {\rm Im}~z_1>0, |{\rm Im}~z_2|=({\rm Im}~z_1)^2\}$.
$C$ is strictly pseudo concave (tube domain with concave boundary), so
for any plurisubharmonic function $\lambda$ defined on the closure of $U$,
${\rm Sup}_U\lambda = {\rm Sup}_{A\cup B}\lambda$. For $p\in U$, let $\mu_p$
be a Jensen measure for $p$ on $A\cup B$. Next set 
$$\lambda_0=-{\rm Im}~z_1+({\rm Im}~z_1)^2+({\rm Im}~z_2)^2 -{1\over 2}({\rm Re}~z_1)^2
-{1\over 2}({\rm Re}~z_2)^2~,$$
$\lambda_0$ is a strictly plurisubharmonic function, $\lambda (0)=0$, and $\lambda < -{r^2\over 2}$
on $B$.
Therefore at least for $p$ close to 0, one must have $\mu_p (A)>0$. Hence if $\lambda=-\infty$
on $A$, one has for $p\in U$ close to 0, $\lambda (p)=-\infty$, and therefore $\lambda =-\infty$
on any connected set containing 0 on which it is defined.
\medskip\noindent
If one wishes to avoid using Jensen measures, one can argue very simply as follows:
Fix $p\in U$ such that $\lambda_0(p)>-{r^2\over 4}$, if as above $\lambda =-\infty$
on $A$, then for $\epsilon>0$ small enough $\lambda_0+\epsilon \lambda$ clearly
violates the maximum principle that in our case is ${\rm Sup}_U={\rm Sup}_{A\cup B}$,
unless $\lambda (p)=-\infty$.
\bigskip\noindent
The application to almost complex structures in real dimension 4 is immediate. Sharper statements
are immediate too, and smoothness requirements are easy to lower (however not to ${\cal C}^1$,
that would be interesting).
\smallskip\noindent
{\bf Proposition 4.} {\it Let $(M,J)$ be a connected almost complex manifold 
of class ${\cal C}^{3,\alpha}$ ($0<\alpha <1$), of real dimension 4. 
Let $p_0\in M$. Let $N$ be a germ of ${\cal C}^{3,\alpha}$ 
2 dimensional surface, $p_0\in N$. Assume that $J(p_0)TN_{p_0}\neq TN_{p_0}$, where
$TN_{p_0}$ is the tangent space to $N$ at $T_0$. Let $\lambda$ be a $J$-plurisubharmonic
function defined on $M$. If $\lambda =-\infty$ on $N$, then $\lambda \equiv -\infty$.}
\bigskip\noindent
Sketch of proof: It is enough to prove that $\lambda=-\infty$ on some non empty open set.
We can then work in a neighborhood of $p_0$. We can assume that $M\subset {\bf C}^2$, $p_0=0$, 
$J(0)=J_{st}$
(the standard complex structure on ${\bf C}^2$) and, after change of variables, that $N$
coincides with ${\bf R}^2$. This change of variables may drop the regularity of $J$
to ${\cal C}^{2,\alpha}$. Furthermore using dilations, as usual, one can assume that $J$
is arbitrarily close in ${\cal C}^2$ norm to the standard complex structure. In that case,
the function $\lambda_0$ considered above is $J$-plurisubharmonic, due to the characterization
of plurisubharmonicity in terms of a Levi form (see e.g. [I-R] Corollary 1.1). Also, with $U$ as above,
$C$ keeps its strong pseudoconcavity property. So the reasoning used in the complex case applies
immediately if one proves that no non constant $J$-plurisubharmonic function on a neighborhood
of $U$ can attain its maximum on $C$. This follows from the fact that
if $q\in C$ there is a $J$-holomorphic disc $u: {\bf D}\to U$, such that $u(0)=q$
and $u(\zeta )\in U$, if $\zeta \neq 0$ (see [I-R] Proposition 5.1).
\bigskip\noindent
{\bf 3) Real dimension 2n.} Of course non $J$-holomorphic surfaces can be $J$-pluripolar.
Indeed ${\bf R}^2\times \{ 0\}$ is pluripolar in ${\bf C}^3$ since it is included
in the $-\infty$ set of $\log |z_3|$. The result here is that generic $n$-dimensional
manifolds $N$ ($TN\cup J(TN)=TM$) are not $J$-pluripolar. We just indicate how to
adapt the reasoning for ${\bf C}^2$ above to the ${\bf C}^n$ case. Define $U$ by:
$${\rm Im}~z_1>0~,~\sum_{j=2}^n({\rm Im}~z_j)^2<({\rm Im}~z_1)^4~,~\sum_{j=1}^n|z_j|^2<r^2~.$$
If again $C$ is the part of the boundary of $U$ corresponding to
$\sum_{j=2}^n({\rm Im}~z_j)^2={\rm Im}~z_1$, then $C$ is of course not strictly pseudo-concave
but at every point its Levi form has at least one negative eigenvalue, and this is all what is needed.

% ppol6

\bigskip\bigskip

\centerline{\bf Appendix. Three Lemmas in Measure Theory.}
\vskip.5truein
\noindent Recall that we say that a  subset $E\subset {\bf C}$ (on which the variable will
now be denoted by $z$) is {\it fat}
if 0 is a point of Lebesgue density of $E$, i.e.
$ { | E\cap \{ |z|<r\} |\over \pi r^2}\to 1\quad {\rm as}~r\to 0~r$, (where, the Lebesgue
measure of a set $F\subset {\bf C}$ is denoted by $|F|$).
\bigskip\noindent
{\bf Lemma A1.} { \it Let $\nu$ be a finite positive measure on ${\bf C}$.
Assume that $\nu (\{ 0\} )=0$, then there is a fat set $E$
such that $z({1\over |z|}\ast \nu)$ tends to 0 as $z\to 0$, $z\in E$.}
\bigskip\noindent
{\bf Lemma A2.} {\it Let $\psi\in L^1$ (on ${\bf C}$). Then, there is a fat set $E$ such that
$z^2(\psi\ast{1\over z^2})$ tends to 0 as $z\to 0$, $z\in E$.}\hfill\break
The convolution makes sense by the theory of singular integrals.
\bigskip\noindent
{\bf Lemma A3.} {\it Let $\rho\in L^1$, and $\epsilon >0$, then the inequality $\rho> {\epsilon
\over |z|^2 |\log |z||}$ cannot be satisfied everywhere on a set $F$ that has positive density at 0,
i.e. satisfying $| F\cap \{ |z|<r\} |>\delta r^2$ for some fixed $\delta >0$, for all $r>0$ small enough.}

\bigskip\noindent
{\bf Remarks.}

\noindent 1) It is enough to prove Lemma A3 for $\epsilon =1$, replacing $\rho$ by
${\rho\over\epsilon}$. Lemma A3 will be enough for us, but obviously much stronger statements 
can be given.
However, there may not exist any fat set on which
$(|z|^2|\log |z||~\rho )$ is bounded.
An example is the following one.
Take a sequence of integers $(n_j)$ such that 
$\sum_j {1\over \sqrt n_j}<+\infty$. For $2^{-n_j-1}<|z|<2^{-n_j}$,
set $\rho (z)=
{ 1\over |z|^2\sqrt{|\log |z||}}$. Elsewhere set $\rho =0$.
\bigskip

\noindent 2) Simple logic shows that for proving Lemma A1, it is enough to prove that
the set $E_\epsilon=\{ z,~ {1\over |z|}\ast \nu \leq {\epsilon\over |z|}\}$ is a fat set, and
again it is enough to prove this for $\epsilon =1$. Then we can define $E$
as follows. For a sequence $r_j>0$ decreasing fast enough to 0 (in particular we require
$r_{j+1}\leq {1\over 2}r_j$), define $E$
by setting $E\cap \{ r_{j+1}<|z|\leq r_j\} = E_{1\over j} \cap \{ r_{j+1}<|z|\leq r_j\}$.
\medskip\noindent
By the same kind of arguments, for establishing Lemma A2, it is  enough to prove that the set
$E_1=\{ z, |\psi\ast{1\over z^2}|\leq {1\over |z|^2}\}$ is a fat set.
\bigskip\noindent
{\bf Proof of Lemma A1.} As said above, we only need to prove that there is a fat set $E$ such that
$|z({1\over |z|}\ast \nu )|\leq 1$ on $E$. By writing $\nu$ as the sum of a measure
with support near 0 and a measure with no mass near zero, the second one playing no role,
we can assume $\| \nu \|\leq 1$.
One has
$${1\over |z|}\ast \nu (z)=\int\int_{|z-t|\geq 2|z|}~~+~~\int_{|z-t|<2|z|}
{1\over |z-t|}~d\nu (t).$$
Trivially $\int\int_{|z-t|\geq 2|z|}{1\over |z-t|}~d\nu (t) \leq {1\over 2 |z|}$.
Fix $\epsilon >0$, we need to prove that for $R>0$ small enough, the set of $z$ satisfying $|z|\leq R$
and where $ \int_{|z-t|<2|z|}
{1\over |z-t|}~d\nu (t)\geq {1\over 2 |z|}$ has measure at most $\epsilon R^2$.
For $R>0$, let $\chi_R$ be the characteristic function
of the set $|z|\leq R$, in ${\bf C}$.
If $|z|\leq R$,
$$\int\int_{|z-t|< 2|z|}{1\over |z-t|}~d\nu (t)\leq A(z)=
\big( {1\over |.|}\chi_{2R}(.)\big) \ast (\chi_{3R}\nu)~(z).$$ 
The standard estimates for $L^1$ norms of convolutions gives us
$$\|  A\|_{L^1} \leq \| {1\over |.|}\chi_{2R}\|_{L^1}~~\| \chi_{3 R}\nu\|~=
~4\pi R~~\|\chi_{3R}\nu\| \leq {\epsilon\over 2} R~,$$
if $R$ is small enough (since $\nu$ has no point mass at 0). 
Then, obviously on the set $|z|\leq R$, $|A(z)|\geq {1\over 2 |z|}$ can occur only on a set
of measure at most $\epsilon R^2$.
\bigskip\noindent
{\bf Proof of Lemma A2.} 
Set $M=\{z~, |\psi\ast {1\over z^2}|>{1\over |z|^2}\}$. It is enough to prove that,
for any fixed $\delta >0$, for $k\in {\bf N}$ large enough
$|M\cap \{ 2^{-k-1}<|z|<2^{-k}\}|\leq \delta 2^{-2k}$. Similarly to above we can write $\psi$ as the sum of
a function supported near 0, of small $L^1$ norm and a function vanishing near 0 and that
plays no role. So we will be able to assume that the $L^1$ norm of $\psi$ is as small as we will need.
By the weak $L^1$ estimates in the Calderon Zygmund theory, we have for some universal constant $C$,
for all $t>0$:
$$|\{z\in {\bf C}~,~|\psi\ast {1 \over z^2}|>t\}|\leq C\|\psi\|_{L^1}{1\over t}.$$
For $2^{-k-1}\leq |z|\leq 2^{-k}$, $|\psi\ast {1\over z_2}|>{1\over |z|^2}$ implies
$|\psi\ast {1\over z^2}|> 2^{2k}$, so the set of $z$ satisfying  $2^{-k-1}\leq |z|\leq 2^{-k}$
and $|\psi\ast {1\over z_2}|>{1\over |z|^2}$ has measure at most $|\leq C\|\psi\|_{L^1}2^{-2k}$
So it is enough to have
$\| \psi \|_{L^1}\leq {1\over C}\delta$, which we can assume, as said above.
\bigskip\noindent
{\bf Proof of Lemma A3.} 
\bigskip\noindent
Assume that $\rho>{1\over |z|^2|\log |z||}$ on a set $F$ satisfying
$|F\cap \{ |z|<r\}|>\delta r^2$, for $r$ small enough. 
Set $\eta=\sqrt{{\delta\over 2\pi}}$ ($<1$). Then

$$|F\cap \{ \eta r< | z| <r\}|>{\delta\over 2}r^2~.$$
Taking $r=\eta^k$ ($k$ large), we get 
$$\int_{F\cap \{ \eta^{k+1} r< | z| <\eta^k\}}\rho > {\delta\over 2 |\log \eta | }~ {1\over k}~,$$
clearly contradicting integrability of $\rho$, by summation over $k$.

% ppol7

\bigskip\bigskip

\centerline{\bf REFERENCES.}
\vskip1truein
\item{[D-S]} K.\ Diederich, A.\ Sukhov.
Plurisubharmonic functions and almost complex Stein structures.
Mich. Math. J. (to appear)
\smallskip
\item{[E]} F.\ Elkhadhra. J-pluripolar subsets and currents on almost complex manifolds.
To appear in Math. Zeit.
\smallskip
\item{[Ga]} T.\ Gamelin. {\it Uniform algebras.} Prentice Hall  1969.
\smallskip
\item{[G-S]} H.\ Gaussier, A.\ Sukhov. Estimates of the Kobayashi metric
on almost complex manifolds. Bull. Soc. Math. France 133 (2005) no 2, 259-273.
\smallskip
\item {[I-R]} S.\ Ivashkovich, J-P.\ Rosay. Schwarz-type lemmas for solutions
of $\overline \partial$-inequalities and complete hyperbolicity of almost
complex manifolds. Ann. Inst. Fourier 54 (2004), 2387-2435.
\smallskip
\item{[McD-S1]} D.\ McDuff - D.\ Salamon. {\it $J$-holomorphic curves and quantum cohomology}. Univ.
Lect. Series AMS, 6 (1994).
\smallskip
\item{[McD-S2]} D.\ McDuff - D.\ Salamon. {\it $J$-holomorphic curves and symplectic topology}.
AMS Colloquium Pub. 52 (2004).
\smallskip
\item{[N-W]} A.\ Nijenhuis, W.\ Woolf. Some integration problems in
almost complex and complex manifolds. Ann. of Math. {\bf 77} (1963), 424--489.
\smallskip
\item{[Pa]}  N.\ Pali. Fonctions plurisousharmoniques et courants
positifs de type (1,1) sur une vari\'et\'e presque complexe.
 Manuscripta Math. 118 (2005) no 3, 311-337.
\smallskip
\item{[Ro1]} J-P.\ Rosay. J-Holomorphic submanifolds are pluripolar.
Math. Z. 253 (2006), 659-665.
\smallskip
\item{[Ro2]} J-P.\ Rosay. Notes on the Diederich-Sukhov-Tumanov normalization for
almost complex structures. Collect. Math. To appear. 

\noindent Preprint: www.mittag-leffler.se/preprints/0708s/info.php?id=08
\smallskip
\item{[Si]} J.-C. Sikorav. Some properties of holomorphic curves in
almost complex manifolds. In {\it Holomorphic Curves in Symplectic
Geometry}, eds. M. Audin and J. Lafontaine, Birkhauser, 1994, 351--361.
\smallskip
\item{[St]} E.L.\ Stout. {\it Polynomial Convexity}. Progress in Math. 261. Birkh\"auser
2007.
\smallskip
\item{[S-T1]} A.\ Sukhov, A.\ Tumanov. Filling hypersurfaces by discs in almost complex
manifolds of dimension 2. Indiana Univ. Math. J.  {\bf 57} (2008), 509-544.

\smallskip 
\item{[S-T2]} A.\ Sukhov, A.\ Tumanov. Filling real hypersurfaces by pseudoholomorphic discs.
J. Geom. Anal {\bf 18} (2008), 632-649.

\smallskip \vskip.5truein\noindent
Jean-Pierre Rosay
\hfill\break
Department of Mathematics\hfill\break
University of Wisconsin\hfill\break
Madison WI 53706, U.S.A.\hfill\break
jrosay@math.wisc.edu
\end